\documentclass{amsart}
\usepackage{amssymb}
\usepackage{amsfonts}
\usepackage{amsmath}
\usepackage{lastpage}
\usepackage[legalpaper,bookmarks=true,colorlinks=true,linkcolor=blue,citecolor=blue]{hyperref}
\usepackage{graphicx}%
\usepackage{fancyhdr}
\usepackage{color}
\usepackage[mathlines]{lineno}
\usepackage{lscape}
\usepackage{epsfig}
\usepackage{natbib}
\usepackage[]{listings}
\usepackage{geometry}
\usepackage{tgbonum}
\fontfamily{qcr}\selectfont

\newtheorem{theorem}{Theorem}
\theoremstyle{plain}

\newtheorem{fact}{Fact}

\newtheorem{proposition}{Proposition}

\numberwithin{equation}{section}

\begin{document}
\Large
\title[Asymptotic laws for record values in the extreme domain of attraction]{Asymptotic laws for upper and strong record values in the extreme domain of attraction and beyond}

\begin{abstract} Asymptotic laws of records values have usually been investigated as limits in type. In this paper, we use functional representations of the tail of cumulative distribution functions in the extreme value domain of attraction to directly establish asymptotic laws of records value, not necessarily as limits in type. Results beyond the extreme value value domain are provided. Explicit asymptotic laws concerning very usual laws are listed as well. Some of these laws are expected to be used in fitting distribution.

\noindent $^{\dag}$ Gane Samb Lo.\\
LERSTAD, Gaston Berger University, Saint-Louis, S\'en\'egal (main affiliation).\newline
LSTA, Pierre and Marie Curie University, Paris VI, France.\newline
AUST - African University of Sciences and Technology, Abuja, Nigeria\\
gane-samb.lo@edu.ugb.sn, gslo@aust.edu.ng, ganesamblo@ganesamblo.net\\
Permanent address : 1178 Evanston Dr NW T3P 0J9,Calgary, Alberta, Canada.\\

\noindent $^{\dag\dag}$ Mohammad Ahsanullah\\
Department of Management Sciences. Rider University. Lawrenceville, New Jersey, USA\\
Email : ahsan@rider.edu\\


\end{abstract}

\maketitle

\section{Introduction} \label{recEvtBeyond_01}

\noindent Let $X$, $X_1$, $X_2$, $\cdots$ be a sequence independent real-valued randoms, defined on the same probability space $(\Omega, \mathcal{A},\mathbb{P})$, with common cumulative distribution function $F$, which has the lower and upper endpoints, and the generalized inverse function respectively defined by

$$
lep(F)=\inf \{x \in \mathbb{R}, \ F(x)>0\}, \  uep(F)=\sup \{x \in \mathbb{R}, \ F(x)<1\} 
$$

\bigskip \noindent and

$$
F^{-1}(u)=\inf \{x \in \mathbb{R}, \ F(x)\geq u\} \ for \ u \in ]0,1[ \ and \ F^{-1}(0)=F^{-1}(0+).
$$

\bigskip \noindent Finally, let us consider the sequence of strong record values $X^{(1)}=X_1$, $X^{(n)}$, $\cdots$ and the sequence of record times  $U(1)=1$, $U(2)$, $\cdots$\\

\noindent Before beginning an asymptotic theory, we should be sure that we have an infinite sequence $(X^{(n)})_{n\geq 1}$. For a bounded random variable with finite upper bound $uep(F)$ such that $\mathbb{P}(X=uep(F))>0$, we have $(X^{(n)}<uep(F))$ finitely often. This happens for classical integer-valued and bounded random variables as Binomial laws. In such cases, the asymptotic theory is meaningless. But, an interesting question would be the characterization the infinite random sequence $(n_k)_{k\geq 1}$ such that $X_{n_k}=uep(F)$ for all $k\geq 1$.\\

\noindent In all other cases, even if  $uep(F)$ is bounded, the sequence $(X^{(n)})_{n\geq 1}$ is infinite. So, the results of this paper apply to \textit{cdf}'s $F$such that $\mathbb{P}(X=uep(F))=0$. In that context, asymptotic laws have been proposed in the literature by many authors like \cite{tata}, \cite{resnick}, \cite{nevzorov}, etc.,  in relation with Extreme Value Theory, as limits in type in the form

\begin{equation}
(\exists (A_n)_{n\geq 1}\subset \mathbb{R}_{+}\setminus \{0\}), \ \exists \ (B_n)_{n\geq 1}\subset \mathbb{R}, \ \frac{X^{(n)}-B_n}{A_n} \rightsquigarrow Z, \label{search1}
\end{equation} 

\bigskip \noindent where $\rightsquigarrow$ stands for the convergence in distribution and $Z$ is a non-degenerate random variable. The motive beneath this search is the following. If we denote by $M(n)=\max(X_1,\cdots,X_n)$ as the $n$-th maximum for $n\geq 1$, it is clear that we have

\begin{equation}
\forall \ n\geq 1, \ X^{(n)}=M(U(n)). \label{search2}
\end{equation} 

\bigskip \noindent Since for any $F$ in the extremal domain of attraction $\mathcal{D}$, we have that for some $\gamma \in \mathcal{R}$,

\begin{equation}
(\exists (a_n)_{n\geq 1}\subset \mathbb{R}_{+}\setminus \{0\}), \ (\exists (b_n)_{n\geq 1}\subset \mathbb{R}, \ \frac{M(n)-a_n}{b_n} \rightsquigarrow Z_{\gamma}, \label{search3}
\end{equation} 

\noindent where the \textit{cdf} of $Z_{\gamma}$ is the Generalized Extreme Value distribution defined by

$$
G_{\gamma}(x)=\exp(-(1+\gamma x)^{1/\gamma}), \ 1+\gamma x>0, \ G_{0}(x)=\exp(-\exp(-x)) \ for \ x\in \mathbb{R}.
$$

\bigskip \noindent In Extreme value Theory, Formula \eqref{search3} is rephrased as \textit{F is attracted by $G_{\gamma}$} denoted by $F \in D(G_{\gamma})$.\\

\noindent From Formulas \eqref{search2} and \eqref{search3} and from the fact that $U(n)\rightarrow +\infty$ as $n \rightarrow +\infty$, the investigation of the validity of \eqref{search1} was justified enough. The results of the cited authors and others were positive with the stunning result that the \textit{cdf} of $Z$ should be of the form $\Phi(g(x))$, $x\in \mathbb{R}$, where $\Phi$ is the \textit{cdf} of the standard normal law and $g$ satisfies one of three definitions (in which
$c$ is a positive constant)

\begin{eqnarray*}
&&g(x)=x, \ x\in \mathbb{R}.\\
&&g(x)=-\infty 1_{(x<0)} + (c \log x) \ 1_{(x\geq 0)}, x\in \mathbb{R}.\\
&& g(x)= (-c \log -x) \ 1_{(x<0)} + \infty \ 1_{(x>0)}, \ x\in \mathbb{R}. 
\end{eqnarray*}

\bigskip \noindent Instead of using this mathematically appealing approach based on functional equations, an other approach consisting in directly finding the asymptotic laws of $X^{(n)}$, not necessarily in the form of Formula \eqref{search1} is possible and we proceed to it here. That approach is based on representations of 
$F \in \mathcal{D}$ of Karamata and de Haan for example.\\

\noindent Our achievement is the finding the asymptotic laws of the records for all $F \in \mathcal{D}$. First, for $\gamma \neq 0$, outside the frame Formula \eqref{search1}, that is as limits in type, and without any further condition. Secondly,  for $\gamma=0$, within the frame of Formula \eqref{search1}, under a general regularity condition. That regularity condition generally holds for usual \textit{cdf}'s.\\

\noindent We also give general conditions to ensure the asymptotic normality of the records values for $F$ not necessarily in the extremal domain. Finally, we give detailed asymptotic laws of the records of a list of remarkable \textit{cdf}'s with specific coefficients.\\

\noindent In this paper we want short, we use many results from Extreme Value Theory and Records Values Theory. So, for more details, we refer the reader to the books of 
\cite{ahsanullah2001}, \cite{nevzorov}, etc for an easy introduction to records and to those of \cite{galambos}, \cite{dehaan}, \cite{resnick}, \cite{lo2018}, etc. concerning Extreme Value Theory.\\

\noindent To finish this introduction, we recall two important tools of extreme value theory that form the basis of our method. The first is the following proposition. Suppose that $X\geq 0$, that is $F(0)=0$. In that case, we define $Y=\log X$ with \textit{cdf} $G(x)=F(e^x)$, $x \in \mathbb{R}$ and we have\\

\begin{proposition} \label{proplo86} (see \cite{lo86}) We have the following equivalences.\\

\noindent (1) If $\gamma>0$,

$$
F \in D(G_{\gamma}) \ \Leftrightarrow\ (G \in D(G_{0}) \ and \ R(x,G)\rightarrow \gamma \ as \ x \rightarrow uep(G)).
$$

\noindent (2) If $\gamma=0$,

$$
F \in D(G_{0}) \ \Leftrightarrow\ (G \in D(G_{0}) \ and \ R(x,G)\rightarrow 0 \ as \ x \rightarrow uep(G)).
$$

\noindent (3) If $\gamma<0$,

$$
F \in D(G_{\gamma}) \ \Leftrightarrow\ (G \in D(G_{\gamma}). 
$$

\end{proposition}

\noindent In the second place, we recall the following representations of \textit{cdf}'s in the extreme value domain that repeatedly will be used in the sequel.

\begin{proposition} \label{portal.rd} (\cite{karamata} and \cite{dehaan}) We have the following characterizations for the three extremal domains.

\bigskip \noindent (a) $F\in D(H_{\gamma })$, $\gamma >0$, if and only if there exist a
constant $c$ and functions $a(u)$ and $\ell (u)$ of $u\rightarrow $ $u\in
]0,1]$ satisfying

\begin{equation*}
(a(u),\ell (u))\rightarrow (0,0)\text{ as }u\rightarrow 0,
\end{equation*}

\bigskip \noindent such that $F^{-1}$ admits the following representation of Karamata%
\begin{equation}
F^{-1}(1-u)=c(1+a(u))u^{-\gamma }\exp \left(\int_{u}^{1}\frac{\ell (t)}{t}dt\right). \label{portal.rdf}
\end{equation}

\bigskip \noindent (b) $F\in D(H_{\gamma })$, $\gamma <0,$ if and only if $uep(F)<+\infty $ and
there exist a constant $c$ and functions $a(u)$ and $\ell (u)$ of $u\in ]0,1]
$ satisfying
\begin{equation*}
(a(u),\ell (u))\rightarrow (0,0)\text{ as }u\rightarrow 0,
\end{equation*}

\bigskip \noindent such that $F^{-1}$ admit the following representation of Karamata

\begin{equation}
uep(F)-F^{-1}(1-u)=c(1+a(u))u^{-\gamma }\exp \left(\int_{u}^{1}\frac{\ell (t)}{t}
dt\right). \label{portal.rdw}
\end{equation}

\bigskip \noindent (c) $F\in D(H_{0})$ if and only if there exist a constant $d$ and a slowly
varying function $s(u)$ such that 
\begin{equation}
F^{-1}(1-u)=d+s(u)+\int_{u}^{1}\frac{s(t)}{t}dt,0<u<1, \label{portal.rdg}
\end{equation}

\bigskip \noindent and there exist a constant $c$ and functions $a(u)$ and $\ell (u)$ of $u\in ]0,1]$ satisfying
\begin{equation*}
(a(u),\ell (u))\rightarrow (0,0)\text{ as }u\rightarrow 0,
\end{equation*}

\bigskip \noindent such that the function $s(u)$ of $u \in ]0,1[$ admits the representation
\begin{equation}
s(u)=c(1+a(u)) \exp \left(\int_{u}^{1}\frac{\ell (t)}{t}dt\right). \label{portal.rdgs}
\end{equation}

\noindent Moreover, if $F^{-1}(1-u)$ is differentiable for small values of $s$ such
that $r(u)=-u(F^{-1}(1-u))^{\prime }=u\ dF^{-1}(1-u)/du$ is slowly varying at
zero, then \ref{portal.rdg} may be replaced by 

\begin{equation}
F^{-1}(1-u)=d+\int_{u}^{u_{0}}\frac{r(t)}{t}dt,0<u<u_{0}<1, \label{portal.rdgr}
\end{equation}

\noindent which will be called a \textit{reduced de Haan representation} of $F^{-1}.$
\end{proposition}

\noindent The rest of the paper is organized as follows. The results are stated in Section \ref{recEvtBeyond_02}. Examples and Applications are given in
Section \ref{recEvtBeyond_03}. The proofs are stated in Section \ref{recEvtBeyond_04}. The computation related to examples in Section \ref{recEvtBeyond_03} are detailed in the Appendix Section \ref{recEvtBeyond_06}. The paper closed by a conclusion in Section \ref{recEvtBeyond_05}.\\

\section{Results} \label{recEvtBeyond_02}

\noindent Before we state our results, we recall that any $F \in \mathcal{D}$ is associated to a couple of functions $(a(u),b(u))$ of $u \in [0,1]$ as defined in the representations of Proposition \ref{portal.rd} for $F \in D(G_{\gamma})$, $\gamma\neq 0$. In the special case where $\gamma=0$, the pair of functions $(a(circ),b(\circ))$ is used in the representation of the function $s(u)$ is f $u \in [0,1]$ in Representation \eqref{portal.rdg}.\\

\noindent We will need the following condition. Let us define for any $n\geq 1$ a finite sum of $n$ standard exponential random variables

$$
S_{(n)}=E_{1,n} + \cdots + E_{n,n},
$$

\noindent denote

$$
V_n=\exp(-S_{(n)}) \ and \ v_n=\exp(-n), \ n\geq 1
$$

\noindent and finally set the hypotheses

$$
(Ha) \ : \ \sup\left\{ \left|\frac{u}{v} - 1 \right|, \ \min(v_n,V_n) \leq u,v \leq \max(v_n,V_n)\right\} \rightarrow_{\mathbb{P}} 0 \ as \ n\rightarrow +\infty,
$$

\noindent

$$
(Hb) \ : \ (\exists \alpha>0), \ \sqrt{n} \ s(v_n) \rightarrow \alpha \ as \ n\rightarrow +\infty,
$$

\noindent \noindent where $\rightarrow_{\mathbb{P}}$ stands for the convergence in probability.\\

\noindent Here are our results that cover the whole extreme value domain of attraction. For $\gamma\neq 0$, we need no condition.\\

\noindent Let us begin by asymptotic laws for $F \in \mathcal{D}$.

\begin{theorem} \label{theoALL1} Let $F \in D(G_{\gamma})$, $\gamma \in \mathbb{R}$. We have :\\

\noindent (a) If $\gamma>0$, the asymptotic law of $X^{(n)}$ is lognormal, precisely

$$
\left(\frac{X^{(n)}}{F^{-1}\left(1-e^{-n}\right)}\right)^{n^{-1/2}} \rightsquigarrow LN(0,\gamma^2),
$$

\noindent where $LN(m,\sigma^2)$ is the lognormal law of parameters $m$ and $\sigma>0$.\\

\noindent (b) If $\gamma>0$ and $X\geq 0$, $Y=\log X \in D(G_0))$ and $R(x,G)\rightarrow \gamma \ as \ x \rightarrow uep(G)$ and we have

$$
\frac{Y^{(n)}- G^{-1}\left(1-e^{-n}\right)}{\sqrt{n}} \rightsquigarrow \mathcal{N}(0,\gamma^2).
$$

\noindent (c) If $\gamma<0$, the asymptotic law of $X^{(n)}$ is lognormal, precisely

$$
\left(\frac{uep(F)-X^{(n)}}{uep(F)-F^{-1}\left(1-e^{-n}\right)}\right)^{n^{-1/2}} \rightsquigarrow \exp(\mathcal{N}(0,\gamma^2)).
$$

\bigskip \noindent (d) Suppose that $\gamma=0$ and  $R(x,G)\rightarrow 0 \ as \ x \rightarrow uep(G)$. If (Ha) and (Hb) hold both, we have

$$
X^{(n)}- F^{-1}\left(1-e^{-n}\right) \rightsquigarrow \mathcal{N}(0,\alpha^2).
$$

\bigskip \noindent More precisely, we have :  Given $\gamma=0$,  $R(x,G)\rightarrow 0 \ as \ x \rightarrow uep(G)$ and (Ha), the above asymptotic normality is valid if and only if (Hb) holds.
\end{theorem}

\bigskip \noindent Beyond distributions in $\mathcal{D}$, we may use the delta-method as follows. Drawing lessons from Theorem \ref{theoALL1}, we might be tempted to generalize point (a) by imposing that $F^{-1}$ satisfies, for some coefficient $\gamma$,

$$
\forall \lambda>0, \ F^{-1}(1-\lambda u)/F^{-1}(1-u)= \lambda^{\gamma} (1+o(1)), \ u \in ]0,1[.
$$

\bigskip \noindent But, by Extreme Value Theory, this would imply that $F \in G_{\gamma}$ and nothing new would happen. But trying a generalization from Point (c) would be successful. Let us define the following hypotheses :\\

\noindent (Ga) $F$ is differentiable in some left neighborhood $]x_0, \ uep(F)[$ of $uep(F)$.\\

\noindent (Gb) The function 

$$
s(x)=e^{-x} \biggr[F^{-1}(1-t)^{}\biggr]^{\prime}_{t=e^{-x}}, \ e^{x}<u_0<1, \ for \ some \ u_0 \in ]0,1[
$$ 

\bigskip \noindent decreases to $0$ as $x \rightarrow +\infty$ and is such that : for any sequence $(x_n,y_n)_{n\geq 1}$ such that

$$
\limsup_{n\rightarrow +\infty} |x_n-y_n|/\sqrt{n} <+\infty,
$$

\bigskip \noindent we have, for some $\alpha>0$,

$$
\lim_{n\rightarrow +\infty} \sqrt{n} \ s\left(\exp(\min(x_n,y_n))\right)=\lim_{n\rightarrow +\infty} \sqrt{n} \ s\left(\exp(\max(x_n,y_n))\right)=\alpha.
$$

\bigskip \noindent We have the following generalization.

\begin{theorem} \label{theoALL2} If $F$ satisfies Assumptions (Ga) and (Gb), we have

$$
X^{(n)}- F^{-1}\left(1-e^{-n}\right) \rightsquigarrow \mathcal{N}(0,\alpha^2)
$$
\end{theorem}

\bigskip \noindent \textbf{Comments}. A firm look at the results shows that for any $F \in \mathcal{D}$, we found the direct asymptotic law of $X^{(n)}$ or that of a function of $X^{(n)}$, mainly $\log X^{(n)}$. For example, Point (d) of Theorem \ref{theoALL1} cannot be applied when $X$ follows a lognormal law but can be applied to $\exp(X)$. This leads to the following rule for all any $F \in \mathbb{D}$ :\\

\noindent (e) If $F \in D(G_{\gamma})$, $\gamma\neq 0$, we apply Points (a) or (c) without any further condition.\\

\noindent (f) If $F \in D(G_{0})$ and $\exp(X) \in D(G_{\gamma})$ for some $\gamma>0$, we apply Point (b) without any further condition.\\

\noindent (g) If $F \in D(G_{0})$ and $s(u) \rightarrow 0$ as $u\rightarrow 0$. If (Ha) and (Hb) holds, we conclude by applying Point (d). If not (as it is
for a lognormal law), we search whether $X_1=\exp(X) \in D(G_{\gamma})$ for some $\gamma>0$ or $X_1=\exp(X)$ fulfills (Ha) and (Hb). If yes, we conclude by Point (b) or by Point (d). If not, we consider $X_2=\exp(X_1)$, and we continue until we reach $X_p=\exp(X_{p-1}) \in D(G_{\gamma})$ for some $\gamma>0$ or $X_p=\exp(X_{p-1})$ for some $p\geq 1$.\\

\section{Examples and applications} \label{recEvtBeyond_03}

\noindent Let us begin to explain how to apply the results for $\gamma=0$. Generally, we may find the function $s(u)$ of $u \in ]0,1[$ by from the $\pi$-variation formula

$$
\forall \ \lambda>0, \ \frac{F^{-1}(1-\lambda u)-F^{-1}(1-u)}{s(u)} \rightarrow -\log \lambda \ as \ u \rightarrow 0.
$$

\bigskip \noindent Another method concerns the special case where $F$ is differentiable on left neighborhood of \textsl{uep(F)}. It is proved in \cite{lo86} that if 
$u\left(F^{-1}(1-u)\right)^{\prime}$ is slowly varying at zero, we have for some $u_0 \in ]0,1[$,

$$
s(u)=-u\left(F^{-1}(1-u)\right)^{\prime} \ for \ u \in ]0,u_0[.
$$

\bigskip \noindent Checking hypothesis (Ha) and (Hb) can be done with the function $s(u)$ of $u\in ]0,1[$, found as explained above.\\
 
\noindent Here are some specific examples. The details for each case is given in the Appendix (Section \ref{recEvtBeyond_06}, \pageref{recEvtBeyond_06}). We begin for light tails :\\

\noindent \textbf{I - $F \in D(G_0)$}.\\

\bigskip \noindent \textbf{(1) $X$ follows an exponential law $\mathcal{E}(\lambda)$, $\lambda>0$}. By Point (b) of Theorem \ref{theoALL1},\\

$$
\frac{X^{(n)}-n}{\sqrt{n}}   \rightsquigarrow \mathcal{N}(0, \lambda^{-2}).
$$

\bigskip \noindent \textbf{(2) $X$ follows a standard normal law $\mathcal{N}(0,1)$}. By Point (d) of Theorem \ref{theoALL1},\\

$$
X^{(n)}-(2 n)^{1/2} \rightsquigarrow \mathcal{N}(0,1/2).
$$

\bigskip \noindent \textbf{(3) $X$ follows a Rayleigh law of parameter $\rho>0$}, with \textit{cdf}

$$
1-F(x)=\exp(-\rho x^2), \ x \geq 0.
$$

\bigskip \noindent By Point (d) of Theorem \ref{theoALL1}, we have\\

$$
X^{(n)}-\left(\frac{n}{\rho}\right)^{1/2}  \rightsquigarrow \mathcal{N}(0,\rho^{-1}/4).
$$

\bigskip \noindent \textbf{(4) $X$ follows the logistic law}, with \textit{cdf}

$$
F(x)=\frac{1}{1+e^{-x}}, \ x \in \mathbb{R}.
$$

\bigskip \noindent By Point (b) of Theorem \ref{theoALL1}, we have\\

$$
\frac{X^{(n)}-n }{\sqrt{n}}  \rightsquigarrow \mathcal{N}(0,1).
$$

\bigskip \noindent \textbf{(5) $X>0$ follows a standard lognormal law, that is $\log X$ follows a standard normal law}. We have\\

$$
\log X^{(n)}-(2n)^{1/2}  \rightsquigarrow \mathcal{N}(0,1/2).
$$

\bigskip \noindent \textbf{(6) $X>0$ a follows a Gumbel law} with \textit{cdf}

$$
F(x)=\exp\left(- e^{-x}\right), \ x\in \mathbb{R}.
$$

\bigskip \noindent By Point (b) of Theorem \ref{theoALL1}, we have\\

$$
\frac{X^{(n)}-n}{\sqrt{n}} \rightsquigarrow \mathcal{N}(0,1).
$$

\bigskip \noindent \textbf{II - $F \in D(G_\gamma)$, $\gamma>0$}.\\

\noindent \textbf{(7) $X$ follows a log-logistic law of parameter $p>0$}, with \textit{cfd}

$$
F(x)=\frac{x^p}{1+x^p}, \ x\geq 0.
$$

\bigskip \noindent By Point (a) of Theorem \ref{theoALL1},\\

$$
\left(e^{-n/p}X^{(n)}\right)^{-1/2}   \rightsquigarrow LN(0, p^{2}).
$$

\bigskip \noindent \textbf{(8) $X$ follows a sing-Maddala law of parameters $a>0$, $b>0$ and $c>0$}, with \textit{cdf}

$$
1-F(x)=\left(\frac{1}{1+ax^b}\right)^{c}, \ x\geq 0.
$$

\bigskip \noindent By Point (a), we have

$$
\left(a^{1/b} \exp(-n/(bc)) X^{(n)}\right)^{1/\sqrt{n}} \rightsquigarrow LN(0, (bc)^{-2}).
$$

\section{Proofs} \label{recEvtBeyond_04}

\noindent \textbf{(I) - Proof of Theorem \ref{theoALL1}}.\\

\noindent We begin by describing the main tools which are based on following results of Records theory. Suppose that $\{T, T_j>0, \ 1\leq j \leq k\}$ are $(k+1)$ non-negative real-valued and \textsl{iid} random variables and define

$$
X_0=0, \ T_j=X_{j}-X_{j-1}, \ 1\leq j \leq k.
$$

\bigskip \noindent It is clear that if $T \sim \mathcal{E}(\lambda)$, $\lambda>0$, then the absolutely continuous \textsl{pdf} of $T=(T_1,\cdots,T_k)^t$ is given by

\begin{equation}
f_T(t_1,\cdots,t_k)=\lambda^k e^{-\lambda t_k} \ 1_{(0\leq t_1\leq \cdots \leq t_k)}. \label{RepExpo}
\end{equation}

\noindent Suppose if $T_j$'s are independent and follow an exponential law $\mathcal{E}(\lambda)$, $\lambda>0$, we have 

$$
r(x)=\frac{dF(x)/dx}{1-F(x)}=\lambda \ and \ f(x)=\lambda e^{-\lambda x}, \ x\geq 0.
$$

\bigskip \noindent As stated in page 3 in \cite{ahsanullah2001}, the joint distribution of the $k$ first records values $(T^{(1)}, \cdots, T^{(k)})$ of the sequence $(T_n)_{n\geq 1}$  is the one given in Formula \eqref{RepExpo}. As a consequence, we have

\begin{fact} \label{factRepExpo} If the the $T_j$'s are independent and follow an exponential law $\mathcal{E}(\lambda)$, the $k$-th record value, $k\geq 1$, has the same law as the sum of $k$ independent $\mathcal{E}(\lambda)$-random variables $E_{1,k}$, $\cdots$, $E_{k,k}$, i.e.

$$
T^{(k)}=_{d} E_{1,k} + \cdots + E_{k,k},
$$
\end{fact}

\bigskip \noindent where $=_d$ stands for the equality in distribution. By the Renyi's representation, we can represented the random variable $X$ of \textit{cdf} $F$ by a standard exponential random variable $E$

$$
X=_{d} F^{-1}\left(1-e^{-E}\right).
$$

\bigskip \noindent It comes that, by considering \textit{iid} sequence  $(X_n)_{n\geq 1}$ and $(E_n)_{n\geq 1}$ from $X$ and $E$ and by denoting the two $n$-th records valued $X^{(n)}$ and $E^{(n)}$ from the two sequences respectively, we have the following representations

$$
X^{(n)}=_{d}F^{-1}\left(1-e^{-E^{(n)}}\right),
$$

\bigskip \noindent where

$$
S_{(n)}=E^{(n)}=E_{1,n} + \cdots + E_{n,n}.
$$

\bigskip \noindent In the sequel, we can and do use the equality : $X^{(n)}=F^{-1}\left(1-e^{-S_{(n)}}\right)$. Let us apply the representations by using the simple central limit theorem

$$
\frac{S_{(n)}-n}{\sqrt{n}} \rightsquigarrow \mathcal{N}(0,1) \ as \ n\rightarrow +\infty.
$$

\bigskip \noindent In the sequel, any unspecified limit is meant as $n\rightarrow +\infty$.\\

\noindent Let us suppose $X \in D(G_{1/\gamma})$. If $X\geq 0$, we will consider $Y = \log X$ of \textit{cdf} $G$ defined by $G(x)=F(e^x)$, $x \in \mathbb{R}$. Let us prove the theorem.\\

\noindent \textbf{(a) - Asymptotic law of $X^{(n)}$ for $\gamma>0$}. We recall that  $V_n=e^{-S_{(n)}}$ and $v_n=e^{-n}$, $n\geq 1$. By Representation  \eqref{portal.rdf}, we have

\begin{equation*}
F^{-1}\left(1-e^{-S_{(n)}}\right)=(1+a(V_n)) V_n^{-\gamma} \exp\left(\int_{V_n}^{1} \frac{b(t)}{t} \ dt\right), \ n\geq 1
\end{equation*}

\noindent and

\begin{equation*}
F^{-1}\left(1-e^{-n}\right)=(1+a(v_n)) v_n^{-\gamma} \exp\left(\int_{v_n}^{1} \frac{b(t)}{t} \ dt\right), \ n\geq 1.
\end{equation*}

\bigskip \noindent We get that $V_n \rightarrow_{\mathbb{P}} 0$, $(1+a(V_n))/(1+a(v_n))\equiv 1+p_n\rightarrow_{\mathbb{P}} 1$. We get

\begin{equation*}
\log\left(\frac{X^{(n)}}{F^{-1}\left(1-e^{-n}\right)}\right)=p_n (1+o_{\mathbb{P}}(1))-\gamma(S_{(n)}-n)+\int_{v_n}^{V_n} \frac{b(t)}{t} \ dt.
\end{equation*}

\bigskip \noindent We have

$$
\left| \int_{v_n}^{1} \frac{b(t)}{t} \ dt\right| \leq \biggr(\sup_{\ 0\leq t \leq (v_n \vee V_n) } |b(t)|\biggr)  |S_{(n)}-n|. 
$$ 

\bigskip \noindent By combining the two later formulae, we have

$$
n^{-1/2}\log\left(\frac{X^{(n)}}{F^{-1}\left(1-e^{-n}\right)}\right) \rightsquigarrow \mathcal{N}(0,\gamma^2).
$$

\bigskip \noindent \textbf{(b) - Asymptotic law of $Y^{(n)}$ for $\gamma>0$}. From the previous theorem, it is immediate for the following result. It is clear the $G^{-1}=\log F^{-1}$. So, the previous theorem implies

\begin{equation*}
n^{-1/2}\left( Y^{(n)} - G^{-1}\left(1-e^{-n}\right)\right) \rightsquigarrow \mathcal{N}(0,\gamma^2).
\end{equation*}

\bigskip \noindent Here, it is clear that $Y \in D(G_0)$ and $R(x,G)\rightarrow \gamma$ as $x \rightarrow uep(G)$. Hence this result says that  

\begin{equation*}
n^{-1/2}\left( X^{(n)} - F^{-1}\left(1-e^{-n}\right)\right) \rightsquigarrow \mathcal{N}(0,\gamma^2),
\end{equation*}

\noindent if $F \in D(G_0)$ and $R(x,F)\rightarrow \gamma$ as $x \rightarrow uep(F)$.\\

\noindent \textbf{(c) - Asymptotic law of $Y^{(n)}$ for $\gamma<0$}. We have $\mathbb{P}(X=uep(F))=0$. By using Representation \eqref{portal.rdw}, we may and do prove this point exactly as for Point (a).\\

\noindent \textbf{(d) - Asymptotic law of $Y^{(n)}$ for $\gamma=0$}. We did not have yet the general law. Let us learn for a no-trivial example.\\

\noindent \textbf{(A) - $X \sim \mathcal{N}(0,1)$}. Let us recall the expansion of the tail of $F$ as follows

\begin{eqnarray}
F^{-1}(1-s) &=&(2\log (1/s))^{1/2}-\frac{\log 4\pi +\log \log (1/s)}{2(2\log (1/s))^{1/2}}  \label{portal.Gauss.quantile} \\
&+&O((\log \log (1/s)^{2}(\log 1/s)^{-1/2})). \notag
\end{eqnarray}

\noindent We have

\begin{eqnarray*}
X^{(n)}&=&(2S_{(n)})^{1/2}-\frac{\log 4\pi +\log S_{(n)}}{2(2S_{(n)}))^{1/2}}+O(S_{(n)}^{-1/2}\log S_{(n)})\\
&=&(2S_{(n)})^{1/2}\biggr(1-\frac{\log 4\pi +\log \log (1/s)}{4S_{(n)}}+ O(S_{(n)}(\log S_{(n)})\biggr)\\
&=&(2S_{(n)})^{1/2}(1+\varepsilon_n).
\end{eqnarray*}

\bigskip \noindent We have that $\varepsilon_n=O_{\mathbb{P}}(n^{-1})$. Let us use the mean value theorem to get 

$$
S_{(n)}^{1/2}-n^{1/2}=\frac{1}{2}\frac{S_{(n)}-n}{\sqrt{n}} (n/\zeta_n)^{1/2},
$$

\bigskip \noindent with $n \wedge S_{(n)} < \zeta_n < n \vee S_{(n)}$ and next, by the weak law of large numbers, $2(S_{(n)}^{1/2}-n^{1/2}) \rightsquigarrow \mathcal{N}(0,1)$. By plugging this in the later formula, we get

\begin{eqnarray}
\biggr(X^{(n)}-(2n)^{1/2}\bigg)-\biggr(\sqrt{2}(S_{(n)})^{1/2}-n^{1/2})\biggr)=O_{\mathbb{P}}(\varepsilon_n) + O_{\mathbb{P}}(n^{-1/2}). \label{normalF}
\end{eqnarray}

\bigskip \noindent We conclude that

$$
X^{(n)}-(2n)^{1/2} \rightsquigarrow \mathcal{N}(0,1/2).
$$

\noindent By putting

$$
b_n=(2 n)^{1/2}-\frac{\log 4\pi +\log n}{2(2n)^{1/2}}
$$

\noindent we also have

$$
X^{(n)}-\biggr((2 n)^{1/2}-\frac{\log 4\pi +\log n}{2(2n)^{1/2}}\biggr) \rightsquigarrow \mathcal{N}(0,1/2)
$$

\bigskip \noindent \textbf{(B) - General proof}. It known that $s(u)\sim R(F^{-1}(1-u),F)$ and so, $s(u)\rightarrow 0$ as $u\rightarrow 0$, By representation \eqref{portal.rdg} of Proposition \ref{portal.rd} and Hypothesis (Ha) together lead to

\begin{eqnarray*}
X^{(n)}-F^{-1}(1-e^n)&=&s(V_n)-s(v_n)+\int_{v_n}^{V_n} \frac{s(u)}{u} \ du\\
&=&o_{\mathbb{P}}(1)+s(V_n)-s(v_n) -(1+o_{\mathbb{P}}(1)) s(v_n) (S_{(n)}-n).
\end{eqnarray*}

\bigskip \noindent From there, the conclusion is immediate. $\blacksquare$\\

\noindent \textbf{(II) - Proof of Theorem \ref{theoALL2}}. We have $g(x)=F^{-1}(1-e^{-x})$ $g^{\prime}(x)=S(x)$, $x \in ]lep(F),uep(F)[$. The mean value theorem gives, for

\begin{equation}
X^{(n)}-F^{-1}(1-e^{-n})=\frac{S_{(n)}-n}{\sqrt{n}} \biggr(\sqrt{n}S(\exp(-\zeta_n))\biggr),
\end{equation}

\bigskip \noindent where

$$
\zeta_n \in ]\min(n,S_{(n)}), \ \max(n,S_{(n)})[.
$$

\bigskip \noindent From there, the conclusion is direct. $\blacksquare$\\

\section{Conclusion} \label{recEvtBeyond_05}

\noindent After the statements of the asymptotic laws of the strong record values from \textit{iid} random variables and after some examples have been given,
it should be interesting to a review of such asymptotic laws for as much as possible \textit{cdf}'s $F \in \mathcal{D}$.  

\newpage
\section{Appendix} \label{recEvtBeyond_06}

\noindent Let us give the details concerning the results listed in Section \ref{recEvtBeyond_03}.\\

\bigskip \noindent \textbf{(1) $X$ follows an exponential law $\mathcal{E}(\lambda)$, $\lambda>0$}. We have $exp(X) \in D(G_{-1})$ and 
$F^{-1}(1-e^{-n})=n$. We apply Point (b) to conclude.\\ 

$$
\frac{X^{(n)}-n}{\sqrt{n}}   \rightsquigarrow \mathcal{N}(0, \lambda^-2).
$$

\bigskip \noindent \textbf{(2) $X$ follows a standard normal law $\mathcal{N}(0,1)$}. The result of this point is justified by Formula \ref{normalF}, page \pageref{normalF}.

\bigskip \noindent \textbf{(3) $X$ a follows Rayleigh law of parameter $\rho>0$}. We have

$$
F^{-1}(1-u)=\left(-\frac{1}{\rho} \log u\right)^{1/2}, \ u \in ]0,1[
$$

\noindent and 

$$
s(u)=-u \left(F^{-1}(1-u)\right)^{\prime}=\frac{1}{2\rho (-(1/\rho) log u)^{1/2}} \rightarrow 0 \ as \ u \rightarrow 0. 
$$ 

\noindent Furthermore, $s(u)$ is decreasing in $u \in ]0,1[$ and $s(V_n)/s(v_n) \rightarrow 0$ as $n\rightarrow +\infty$. Finally,

$$
\sqrt{n} s(v_n) \rightarrow \rho^{-1/2}/2.
$$

\bigskip \noindent We conclude the case by applying Point (d) of Theorem \ref{theoALL1}.\\

\bigskip \noindent \textbf{(4) $X$ a follows the logistic law}. It is immediate that $\exp(X) \in D(G_{-1})$ and we have

$$
F^{-1}(1-u)=\log(u/(1-u)), \ u \in ]0,1[.
$$

\noindent We conclude with Point (b) of Theorem \ref{theoALL1}.\\

\bigskip \noindent \textbf{(5) $X>0$ a follows a standard lognormal law, that is $\log X$ follows a standard normal law}.\\

\noindent Since $\log X^{(n)}$ has the same law as the $n$-th record $Z^{(n)}$ from \textit{iid} $\mathcal{N}(0,1)$ random variables. So we have

$$
\log X^{(n)} - (2n)^{1/2} \rightarrow \mathcal{0}(0,1/2).
$$

\bigskip \noindent \textbf{(6) $X>0$ a follows a Gumbel law}. We have

$$
F^{-1}(1-u)=-\log \log(1/(1-u)), \ u \in ]0,1[
$$

\bigskip \noindent and for any $\lambda>0$.

\begin{eqnarray*}
F^{-1}(1-\lambda u)-F^{-1}(1-u)=\log(\lambda(1+o(1))) \rightarrow \log \lambda \ as \ u \rightarrow 0.
\end{eqnarray*}

\bigskip \noindent So, $\exp(X) \in D(G_{-1})$. From there, an application of Point (b) of Theorem \ref{theoALL1} closes the case.\\

\bigskip \noindent \textbf{(7) $X$ follows a log-logistic law of parameter $p>0$}, with \textit{cfd}

$$
F(x)=\frac{x^p}{1+x^p}, \ x\geq 0.
$$

\noindent We have

$$
F^{-1}(1-u)= u^{-1/p} (1-u)^{1/p}, u \in ]0,1[.
$$

\bigskip \noindent By Point (a) of Theorem \ref{theoALL1},\\

$$
\left(e^{-n/p} X^{(n)}\right)^{1/\sqrt{n}}   \rightsquigarrow LN(0, p^{2}).
$$

\bigskip \noindent \textbf{(8) $X$ follows a sing-Maddala law of parameters $a>0$, $b>0$ and $c>0$}. We have

$$
1-F(x)=x^{-bc} (x^{-b}+a)^{-c} \equiv x^{-bc} L(x), \ x \geq 0,
$$

\noindent and  $L$ is a slowly varying function at $+\infty$. So $F \in G_{1/(bc)}$. Applying of Point (a) of Theorem, when combined with

$$
F^{-1}(1-u)=a^{-1/b} u^{-1/(bc)}(1-u^{1/c})^{1/b}, \ u \in ]0,1[,
$$

\noindent and with,

$$
F^{-1}(1-e^{-n})=a^{-1/b} e^{n/(bc)}(1-e^{-n/c})^{1/b},
$$

\bigskip \noindent for $n\geq 1$, closes the case.\\

\end{document}